\documentclass{article}

\usepackage{amssymb}
\usepackage{amsmath}
\usepackage[dvips]{graphicx}

\def\be{\begin{equation}} 
\def\ee{\end{equation}}
\newcommand{\beq}{\begin{eqnarray}}
\newcommand{\eeq}{\end{eqnarray}}
\newcommand{\nbeq}{\begin{eqnarray*}}
\newcommand{\neeq}{\end{eqnarray*}}
\def\D{\displaystyle}

\begin{document}
\title{ Characterizations of distributions via order statistics with random exponential shifts}
\author{ M. Ahsanullah \and  V.B. Nevzorov  \and G.P. Yanev 
}
\date{\empty}
\maketitle

\begin{abstract}
We characterize probability distributions via equalities in law between two order statistics shifted by independent exponential variables. An explicit formula for the quintile function of the identified family of distributions is obtained. The results extend some known characterizations of exponential and logistic distributions.

\vspace{0.5cm}\noindent {\bf Keywords} characterizations, order statistics, exponential distribution, generalized logistic distribution.
\end{abstract}

\section{Introduction}

Distributional (in law) relations between order statistics are useful and elegant tools for characterizing probability distributions. More specifically, there is a large number of publications on characterizations based on recurrences involving both order statistics and standard exponential variables. An excellent review and discussion of the available results on this subject can be found in \cite{GKS98}.

AlZaid and Ahsanullah \cite{AA03}  showed that a non-negative random variable $X$ with an absolutely continuous cumulative distribution function (cdf) has standard exponential distribution iff for a fixed $k$, such that  $1\le k\le n-1$,
\be \label{AA03}
X_{k,n}+\frac{\xi}{n-k}\stackrel{d}{=}X_{k+1,n},
\ee
where $X_{1,n}, X_{2,n}\ldots , X_{n,n}$ are the order statistics in a sample with parent $X$ and $\xi$ is standard exponential and independent from $X_{k,n}$. Please, see also \cite{WA04} and \cite{NA08} for an alternative proof and some extensions in the context of random contractions (cf. \cite{N01}, p.14). On the other hand, it was shown in Ahsanullah et al. \cite{AYO11}, under some additional regularity assumptions, that a random variable $X$ with absolutely continuous cdf is standard logistic if and only if for a fixed $k$ such that $1\le k\le n-1$,
\be \label{AYO11}
X_{k,n}+\frac{\xi_1}{n-k}\stackrel{d}{=}X_{k+1,n}-\frac{\xi_2}{k},
\ee
where
$\xi_{1}$ and $\xi_{2}$ are independent standard exponential variables, which are also independent from $X_{k,n}$ and $X_{k+1,n}$, respectively.

First, we will study some extensions of the distributional relations (\ref{AA03}) and (\ref{AYO11}) involving order statistics shifted by independent exponential random variables. What can be said about the distribution of the parent variable $X$ if the following more general distributional equality holds
\be \label{order stats gen}
X_{k,n}+a\xi_1 \stackrel{d}{=} X_{r,n}-b\xi_2, \qquad  a\ge 0,\quad b\ge 0, \quad 1\le k<r\le n,
\ee
where $\xi_{1}$ and $\xi_{2}$ are independent standard exponential variables, which are also independent from $X_{k,n}$ and $X_{r,n}$, respectively?

Our first result answers the above question in the case of adjacent order statistics, i.e., when $r=k+1$. Throughout this paper we will use the term $Q(y)$ to refer to the pseudo-inverse (quintile) function of $F(x)$, i.e., $Q(y)=\inf\{x:\ F(x)\ge y\}$ for  $y\in(0,1)$.

\vspace{0.3cm}{\bf Theorem 1}\ {\it Let $k$ be a fixed integer such that $1\le k\le n-1$ and  $a\ge 0$ and $b\ge 0$ be two fixed real numbers such that $(a,b)\ne (0,0)$. Assume $X$ is a random variable with  continuous cdf $F(x)$ and $\xi_{1}$ and $\xi_{2}$ are independent standard exponential variables, which are also independent from $X_{k,n}$ and $X_{k+1,n}$, respectively.
Then
\be \label{order stats}
X_{k,n}+a\xi_1\stackrel{d}{=} X_{k+1,n}-b\xi_2
\ee
holds 
if and only if the quintile function $Q(x)$ satisfies
\be \label{mainQ}
Q(x)=\log\left\{ cx^{ bk}(1-x)^{ -a(n-k)}\right\}, \qquad 0<x<1,
  \ee
where $c>0$ is an arbitrary constant.
}

{\bf Remarks.}\ The following particular cases of (\ref{order stats}) and (\ref{mainQ}) might be of independent interest. (i) Setting $a=1/(n-k)$ and $b=0$ in (\ref{mainQ}), we obtain $F(x)=1-e^{-(x-c)}$ for $x\ge c$, i.e., an exponential cdf.
(ii) If $a=0$ and $b=1$ we have $F(x)=\exp\{x-c\}$ for $-\infty<x\le c$.
(iii) If $bk=a(n-k)$, then (\ref{mainQ}) yields the logistic cdf (cf. \cite{JKB95}, p.114)
\[
F(x)=\frac{c_1}{c_1+\exp\{-x\}}, \qquad -\infty<x<\infty,\quad c_1>0.
\]
(iv) If $n=2k$ and $a=b\ne 0$, then Theorem 1 implies
\[
X_{k,2k}+a\xi_1\stackrel{d}{=} X_{k+1,2k}-a\xi_2\quad \mbox{iff}\quad F(x)=\frac{1}{1+c_2\exp\{-x/k\}},\quad c_2>0.
\]

\vspace{0.3cm}The next result addresses the case when (\ref{order stats gen}) involves non-adjacent (two spacings away) order statistics.

\vspace{0.3cm}{\bf Theorem 2}\ {\it  Let $k$ be a fixed integer such that $1\le k\le n-2$ and  $a\ge 0$ and $b\ge 0$ be two fixed real numbers such that $(a,b)\ne (0,0)$. Assume $X$ is a random variable with  continuous cdf $F(x)$ and $\xi_{1}$ and $\xi_{2}$ are independent standard exponential variables, which are also independent from $X_{k,n}$ and $X_{k+2,n}$, respectively. Then
\be \label{order stats2}
X_{k,n}+a\xi_1 \stackrel{d}{=}X_{k+2,n}-b\xi_2
\ee
holds 
if and only if the quintile function $Q(x)$ satisfies
\be \label{mainQ21}
Q(x)=\log\left\{ cx^{bk}(1-x)^{-a(n-k-1)}\right\}-W(x; k,n), \qquad 0<x<1,
  \ee
where  $c>0$ is an arbitrary constant and
\be \label{W}
W(x; k,n)=\left\{
  \begin{array}{ll}
   \frac{\D d\log\{(n-2k+1)x+k+1\}}{\D n-2k+1}, & \mbox{if}\quad n\ne 2k+1; \\
\\
   (a+b)kx, & \mbox{if}\quad n=2k+1,
  \end{array}
\right.
\ee
where $d=bk(n-k)+a(n-k-1)(k+1)$.
}

{\bf Remark.}\ Note that if $n=2k+1$ and $a=b=1/k$, then Theorem 2 implies
\[
X_{k,2k+1}+\frac{\xi_1}{k}\stackrel{d}{=}X_{k+2,2k+1}-\frac{\xi_2}{k}\quad \mbox{iff}\quad Q(x)=\log\left\{\frac{cx}{1-x}\right\}-2x.
\]

A number of characterizations of the logistic distribution use distributional relations  between $X$ and order statistics with positive and negative exponential random shifts. George and Mudholkar \cite{GM81} (see also \cite{LH08}) proved that $X$ is standard logistic if and only if
\[
X\stackrel{d}{=} X_{1,2}+\xi  \qquad \mbox{or}\qquad X\stackrel{d}{=}X_{2,2}-\xi.
\]
More generally, Lin and Hu \cite{LH08} established that, under some smoothness conditions,  $X$ is standard logistic if and only if
\[
X\stackrel{d}{=} X_{1,n}+\sum_{j=1}^{n-1} \frac{\xi_j}{j},
\]
where $\xi_j$ for $j=1,2,\ldots, n-1$ are independent standard exponential variables, which are also independent from $X_{1,n}$. Ahsanullah et al. \cite{AYO11} proved that the standard logistic distribution is characterized by the following distributional equality holding for a fixed $k$, such that $1\le k\le n-1$,
\[
X\stackrel{d}{=} X_{1,n}+\sum_{j=1}^{n-k} \frac{\xi_1(j)}{j}-\sum_{j=1}^{k-1} \frac{\xi_2(j)}{j},
\]
where $X_{1,n}$, $\xi_1(j)$, and $\xi_2(j)$ for $j=1,2,\ldots ,n-1$ are mutually independent and all $\xi$'s are standard exponential. Finally, Zykov and Nevzorov \cite{ZN10} obtained characterizations based on either
\be \label{ZN10}
X\stackrel{d}{=}X_{n,n}-\xi\qquad  \mbox{or} \qquad  X + \xi \stackrel{d}{=}X_{n,n}.
\ee

\vspace{0.3cm}Exploring the so-called $F^{\alpha}$-scheme (cf. \cite{N01}, Lecture 25), we will study characterizations of logistic and related distributions based on one extension of (\ref{ZN10}).

\vspace{0.3cm} {\bf Theorem 3}  {\it  Let $a\ge 0$ and $b\ge 0$ be two fixed real numbers such that $(a,b)\ne (0,0)$ and $F(x)$ be a  continuous distribution function. Suppose $Y_{1}$ and $Y_{2}$ are independent random variables with distribution functions $F^\alpha(x)$ and $F^\beta(x)$, respectively for  $\alpha >0$ and  $\beta >0$. Furthermore, let $\xi_{1}$ and $\xi_{2}$ be independent standard exponential variables, which are also independent from $Y_1$ and $\max\{Y_1,Y_2\}$, respectively.
Then
\be \label{main}
Y_1+a\xi_1\stackrel{d}{=} \max\{Y_1,Y_2\}-b\xi_2
\ee
holds 
if and only if the quintile function $Q(x)$ satisfies
\be \label{new}
Q(x)=\log\left\{ cx^{ b\alpha }(1-x^\beta)^{ -d}\right\}, \qquad 0<x<1,
  \ee
where $c>0$ is an arbitrary constant and $d=[a(\alpha+\beta)+b\alpha]/\beta $.}

{\bf Remarks.}\ (i) Setting $a=1$ and $b=0$ in (\ref{main}), yields the cdf
\[
F(x)=\left(1-\exp\left\{-\frac{\beta}{\alpha+\beta} (x-c)\right\}\right)^{1/\beta}, \qquad  c\le x<\infty.
 \]
(ii) If $b\alpha=\beta$ and $d=1$, then (\ref{new}) is the quintile function of the Type~I generalized logistic cdf (cf. \cite{JKB95}, p. 140)
\[
F(x)=\left(\frac{c_2}{c_2+\exp\{-x\}}\right)^{1/\beta}, \qquad -\infty<x<\infty,\quad c_2>0.
\]

\vspace{0.3cm}The special case of Theorem 3 when  $\alpha=1$ and $\beta=n-1$ is particularly noteworthy.

\vspace{0.3cm}{\bf Corollary}\ {\it  Let $a\ge 0$ and $b\ge 0$ be two real numbers such that $(a,b)\ne (0,0)$. Assume $X$ is a random variable with  continuous cdf $F(x)$ and $\xi_{1}$ and $\xi_{2}$ are independent standard exponential variables, which are also independent from $X$ and $X_{n,n}$, respectively.
Then
\be \label{max}
X+a\xi_1\stackrel{d}{=}X_{n,n}-b\xi_2, \qquad n\ge 2,
\ee
holds 
if and only if the quintile function $Q(x)$ satisfies
\be \label{mainQ2}
Q(x)=\log\left\{ cx^{ b}(1-x^{n-1})^{ -d}\right\}, \qquad 0<x<1,
  \ee
where $c>0$ is an arbitrary constant and $d=(b+na)/(n-1)$.
}

The rest of the paper is organized as follows. In Section 2, we present an auxiliary result before proving Theorem 1. Section 3 deals with the case of non-adjacent order statistics. The proof of Theorem 3 is given in Section 4. We summarize the findings and discuss potential future work  in the last section.

\section{Adjacent order statistics} Let $Z_1$ and $Z_2$ be two random variables with  continuous distribution functions $G$ and $H$, respectively. Let $\xi_1$ and  $\xi_2$ be independent standard exponential variables, which are also independent from $Z_1$ and $Z_2$, respectively.

\vspace{0.3cm}{\bf Lemma}\ {\it Let $a$ and $b$ be two real numbers such that $(a,b)\ne (0,0)$. Then
\be \label{lemma1}
Z_1+a\xi_1 \stackrel{d}{=} Z_2-b\xi_2,
\ee
if and only if
\be \label{lemma2}
bG'(x)+aH'(x)  = G(x)-H(x).
\ee
}

\vspace{0.3cm}{\bf Proof.}\ Let us first prove that (\ref{lemma1}) implies (\ref{lemma2}).
Assume both $a\ne 0$ and $b\ne 0$. It is not difficult to see that (\ref{lemma1}) is equivalent to
\[
\frac{1}{a}\int_{-\infty}^x G(u)\exp\left\{\frac{ u-x}{ a}\right\}du=\frac{1}{b}\int_x^\infty H(u)\exp\left\{\frac{x-u}{b}\right\}du,
\]
which, in turn,  implies
\[
b\int_{-\infty}^xG(u)\exp\left\{\frac{u}{a}\right\}du=a\exp\left\{\frac{(a+b)x}{ab}\right\}\int_x^\infty H(u)\exp\left\{-\frac{u}{b}\right\}du.
\]
Since $G$ and $H$ are continuous, we can differentiate with respect to $x$ both sides of the above equation and  obtain
\[
\exp\left\{-\frac{x}{b}\right\}\left[bG(x)+aH(x)\right]=\frac{a+b}{b}\int_x^\infty H(u)\exp\left\{-\frac{u}{b}\right\}du.
\]
Differentiating with respect to $x$ again (using the  continuity of $H(x)$) we obtain (\ref{lemma2}).
The proof when exactly one of $a$ and $b$ equals zero is similar and is omitted.

To prove that (\ref{lemma2}) implies (\ref{lemma1}), one needs to follow the steps in the above proof in reverse order. The proof of the lemma is complete.

\vspace{0.3cm} Further on we denote  $F_{i,n}(x)=P(X_{i,n}\le x)$ for $1\le i\le n$. Recall that
\be \label{OS_cdf}
F_{i,n}(x)=\sum_{j=i}^n {n \choose j}F^j(x)(1-F(x))^{n-j}, \qquad 1\le i\le n.
\ee

\vspace{0.3cm}{\bf Proof of Theorem 1.}\
First, we will prove that (\ref{order stats}) implies (\ref{mainQ}).
Assume both $a\ne 0$ and $b\ne 0$. It follows from Lemma 1 with $G(x)=F_{k,n}(x)$ and $H(x)=F_{k+1,n}(x)$, making use of (\ref{OS_cdf}), that
\beq \label{eqn1}
bF'_{k,n}(x)+aF'_{k+1,n}(x) & = & F_{k,n}(x)-F_{k+1,n}(x)\\
    & = & {n \choose k}F^k(x)(1-F(x))^{n-k}. \nonumber
\eeq
On the other hand, for the left-hand side of (\ref{eqn1}), using again (\ref{OS_cdf}),  we obtain
\beq \label{eqn2}
\lefteqn{bF'_{k,n}(x)+aF'_{k+1,n}(x)}\\
      & = & (a+b)\frac{d}{dx}\left[\sum_{j=k+1}^{n} {n \choose j}F^j(x)(1-F(x))^{n-j}\right]
 +b\frac{d}{dx}\left[{n \choose k}F^{k}(x)(1-F(x))^{n-k}\right]\nonumber  \\
    & = &
    (a+b)n{n-1 \choose k}F^k(x)(1-F(x))^{n-k-1}F'(x)\nonumber \\
& & +bnF^{k-1}(x)(1-F(x))^{n-k-1}F'(x)\left[{n-1 \choose k-1}(1-F(x))-{n-1 \choose k}F(x)\right]\nonumber \\
    & = &
an{n-1 \choose k}F^k(x)(1-F(x))^{n-k-1}F'(x)+bn{n-1 \choose k-1}F^{k-1}(x)(1-F(x))^{n-k}F'(x)\nonumber \\
    & = &
{n \choose k}F^{k-1}(x)(1-F(x))^{n-k-1}F'(x)[a(n-k)F(x)+bk(1-F(x))].\nonumber
\eeq
Equating the right-hand sides of (\ref{eqn1}) and (\ref{eqn2}), and dividing throughout by  ${n\choose k}F^k(x)(1-F(x))^{n-k}\ne 0$, we obtain
\[
bk\frac{F'(x)}{F(x)}+a(n-k)\frac{F'(x)}{1-F(x)}=1,
\]
which, upon integration with respect to $x$, yields
\be \label{last1}
\log \left\{ F^{bk}(x)(1-F(x))^{-a(n-k)}\right\}=x-c,
\ee
for an arbitrary constant $c>0$. Finally, replacing $x$ with $Q(x)$ in the last equation, we obtain (\ref{mainQ}).
The proof for the cases when either $a$ or $b$ is zero is similar and is omitted here.

To prove the "if" part of the theorem, i.e., (\ref{mainQ}) implies (\ref{order stats}), one needs to start with (\ref{last1}) and repeat the steps of the above proof going backwards to (\ref{eqn1}).

\section{Non-adjacent order statistics}

In the beginning we will prove that (\ref{order stats2}) yields (\ref{mainQ21}). Assume both $a\ne 0$ and $b\ne 0$. Referring to (\ref{OS_cdf}), we have
\beq \label{eqn11}
\lefteqn{\hspace{-0.7cm}F_{k,n}(x)-F_{k+2,n}(x)}\\
    & & =  {n \choose k}F^k(x)(1-F(x))^{n-k} + {n \choose k+1}F^{k+1}(x)(1-F(x))^{n-k-1}. \nonumber
\eeq
On the other hand, similarly to (\ref{eqn2}), we obtain
\beq \label{eqn12}
bF'_{k,n}(x)+aF'_{k+2,n}(x)
  & = &  bk{n \choose k}F^{k-1}(x)(1-F(x))^{n-k}F'(x)\\
&  + & a(n-k-1){n\choose k+1}F^{k+1}(x)(1-F(x))^{n-k-2}F'(x). \nonumber
\eeq
Lemma~1 with $G(x)=F_{k,n}(x)$ and $H(x)=F_{k+2,n}(x)$ implies that the left-hand sides of (\ref{eqn11}) and (\ref{eqn12}) are equal. Equating the right-hand sides of these equations, it is not difficult to see that
\[
F'(x)=\frac{(k+1)F(x)(1-F(x))^2+(n-k)F^2(x)(1-F(x))}{bk(k+1)(1-F(x))^2+a(n-k-1)(n-k)F^2(x)}.
\]
Replacing $x$ with $Q(x)$ and using that $F'(Q(x))=1/Q'(x)$, we obtain
\be \label{Qeqn}
Q'(x)=\frac{a(n-k-1)(n-k)x^2+bk(k+1)(1-x)^2}{x(1-x)(k+1+(n-2k-1)x)}.
\ee
If $n\ne 2k+1$, then (\ref{Qeqn}) can be written as
\be \label{last2}
Q'(x)
    =  \frac{a(n-k-1)}{1-x}+\frac{bk}{x}-\frac{a(n-k-1)(k+1)+bk(n-k)}{k+1+(n-2k-1)x}.
\ee
After integrating with respect to $x$ we obtain (\ref{mainQ21}) with the upper branch of (\ref{W}).
In the case $n= 2k+1$, (\ref{Qeqn}) becomes
\be \label{last3}
Q'(x) 
    =
    \frac{ak}{1-x}+\frac{bk}{x}-(a+b)k
    \ee
and (\ref{mainQ21}) follows after integrating with respect to $x$ again.
The proof when exactly one of $a$ and $ b$ equals zero is omitted.

Starting with (\ref{last2}) and (\ref{last3}) and following the steps of the above proof going backwards, one can verify that (\ref{mainQ21}) yields (\ref{order stats2}).

\section{$F^\alpha$ scheme}

Let us first prove that (\ref{main}) implies (\ref{new}).
Assume $a\ne 0$ and $b\ne 0$. Furthermore, assume $\beta=1$. It follows from Lemma~1 with $G(x)=F^\alpha(x)$ and $H(x)=F^{\alpha+1}(x)$ that
\[
b\alpha F^{\alpha-1}(x)F'(x)+a(\alpha+1)F^{\alpha}(x)F'(x)=F^{\alpha}(x)-F^{\alpha+1}(x),
\]
which, upon division by $F^{\alpha-1}(x)\ne 0$, becomes
\[
F'(x)\left[b\alpha+a(\alpha+1)F(x)\right]=F(x)(1-F(x)).
\]
This, for $F(x)(1-F(x))\ne 0$, can be written as
\[
\frac{1}{F'(x)}=\frac{b\alpha}{F(x)}+\frac{(a+b)\alpha+a}{1-F(x)}.
\]
Replacing $x$ with the quintile function $Q(x)$, we obtain
\[
Q'(x)=\frac{b\alpha}{x}+\frac{(a+b)\alpha+a}{1-x}.
\]
Integrating both sides with respect to $x$, we have
\be \label{beta1}
Q(x)=  \log \left\{cx^{b\alpha}(1-x)^{-(a+b)\alpha-a}\right\},
\ee
where $c>0$ is an arbitrary constant. This completes the proof when $\beta=1$. To extend this proof to the general case $\beta>0$, consider the independent random variables $Y'_1$ and $Y'_2$ with cdf's $(H(x))^{\alpha/\beta}$ and $H(x)$, respectively, where $H(x^\beta)=F(x)$. It follows from (\ref{beta1}) that the quintile function $H^{-1}(x)$ satisfies
\be \label{beta=12}
H^{-1}(x)= \log \left\{cx^{b\alpha/\beta}(1-x)^{-(a+b)\alpha/\beta-a}\right\}.
\ee
Now, taking into account that $Q(x)=H^{-1}(x^\beta)$, we see that (\ref{beta=12}) yields (\ref{new}).
The proof of the cases when either $a$ and $b$ is zero is omitted.

It is not difficult to verify that (\ref{new}) implies (\ref{main}), following the steps of the above proof in reverse order.

\section{Concluding remarks}

We obtained characterizations based on distributional equalities between two order statistics plus or minus  multiples of independent standard exponential variables. The resulting family of distributions includes as its members exponential, logistic, and generalized logistic. In case of adjacent order statistics, the quintile function of the underlying distribution takes on a compact explicit form. If the order statistics are two spacings away, the expression for the quintile function includes an additional additive term, which is linear if the order statistics are on both sides of the median. Further calculations show that if higher-order spacings are involved then the resulting quintile function will have more additive terms.

Exploring the so called $F^\alpha$-scheme, we studied the distribution of a continuous $X$ for which $X$ plus a multiple of a standard exponential variable equals the maximum order statistic minus a multiple of another standard exponential variable. The obtained results generalize those of Zykov and Nevzorov (2010).

One area of future work will be to study in more detail the distributions with quintile functions (\ref{mainQ}), (\ref{mainQ21}), and (\ref{new}). It also remains to be seen what the distribution of $X$ is if the corresponding order statistics differ by a linear combination of standard exponential variables. Some results in this direction are given in Ahsanullah et al. \cite{AYO11}.

\vspace{0.3cm}{\bf Acknowledgements}\ The work of the second author was partially supported by
RFBR grants 10-01-00314a , 09-01-00808 and by FCP grant 2010-1.1 -111-128-033.

\end{document}